\documentclass[11pt,a4paper,twoside]{article}
\usepackage[UKenglish]{babel}
\usepackage[T1]{fontenc}
\usepackage{latexsym,amsfonts,amsmath,amsthm,amssymb,mathrsfs}
\usepackage{fullpage}
\usepackage{paralist}
\usepackage{todonotes}
\usepackage{dsfont}
\usepackage{float}
\usepackage{caption}
\usepackage{xfrac}

\usepackage{graphicx}

\usepackage{fourier}

\newcommand{\yn}{(y_n)_{n \in \mathbb{N}} }
\newcommand{\kj}{n_{\ell}}

\DeclareMathOperator{\R}{\mathbb{R}}
\DeclareMathOperator{\N}{\mathbb{N}}
\DeclareMathOperator{\Q}{\mathbb{Q}}
\DeclareMathOperator{\T}{\mathbb{T}}

\DeclareMathOperator{\Var}{Var}

\newtheorem{thm}{Theorem}

\newtheorem{lem}{Lemma}[section]

\newtheorem{df}[lem]{Definition}
\newtheorem{proposition}[lem]{Proposition}
\newtheorem{rem}[lem]{Remark}

\newtheorem*{rem*}{Remark}

\allowdisplaybreaks
\title{\bf On Birkhoff sums that satisfy no temporal distributional limit theorem for almost every irrational}
\medskip

\author{Lorenz Fr\"uhwirth and Manuel Hauke}

\date{}

\begin{document}

\maketitle
\begin{abstract}
Dolgpoyat and Sarig showed that for any piecewise smooth function $f: \T \to \R$ and
almost every pair $(\alpha,x_0) \in \T \times \T$, $S_N(f,\alpha,x_0) := \sum_{n =1}^{N} f(n\alpha + x_0)$ fails to fulfill a temporal distributional limit theorem. In this article, we show that the two-dimensional average is in fact not needed: For almost every $\alpha \in \T$ and all $x_0 \in \T$, $S_N(f,\alpha,x_0)$ does not satisfy a temporal distributional limit theorem, regardless of centering and scaling. The obtained results additionally lead to progress in a question posed by Dolgopyat and Sarig.
\end{abstract}

\section{Introduction and main results}

Let $X$ be a metric space, $T : X \to X$ a Borel measurable map, $f: X \to \R$ a measurable function and $x_0 \in X$. 
Then

\[S_N(f,T,x_0) = \sum_{k=0}^{N-1} f \circ T^k(x_0)\]
defines the Birkhoff sum of $f$ over $T$ at stage $N$ with starting point $x_0$.
A pair $(T,f)$ is said to satisfy a temporal distributional limit theorem (TDLT) along the orbit of a fixed $x_0 \in X$ whenever there exist two sequences 
$(A_M(f,T,x_0))_{M \in \mathbb{N}}$, $(B_M(f,T,x_0))_{M \in \mathbb{N}}$ with $\lim_{M \to \infty} B_M = \infty$, and a non-constant random variable $Y$ such that

\begin{equation}\label{TLT}\lim_{M \to \infty} \frac{1}{M} \# \left\{1 \leq N \leq M: \frac{S_N(f,T,x_0) - A_M}{B_M} \leq a\right\} = \mathbb{P}[Y \leq a].
\end{equation}
For a more detailed introduction in this area, we refer the reader to \cite{brom_ulc} and especially to the survey article \cite{dol_fay}.\\

Motivated by various research areas such as Discrepancy theory (see, e.g., \cite{beck2,beck3,schmidt}) and the theory of ``deterministic random walks'' (see, e.g., \cite{aake,addo}), particularly interesting and well-studied objects are ergodic sums induced by the irrational rotation on the torus $\T $ (see Section \ref{Sec_preReq} for notation and precise definitions)

\[
\begin{split}
T_{\alpha}: \T & \to \T
\\x &\mapsto x + \alpha,
\end{split}
\]
where $\alpha \notin \Q$.
The corresponding sum $S_N(f,\alpha,x_0) := S_N(f,T_\alpha,x_0)$ is often known as the Birkhoff sum of the irrational circle rotation.\\

There are two different types of temporal limit laws, which we define by following the definition in 
\cite{quenched} as ``quenched'' and ``annealed''.
In the annealed case, the average is not only taken over $N$ for fixed $\alpha$, but a pair
$(\alpha,N)$ is drawn uniformly at random from $\mathbb{T} \times \{1,\ldots,M\}$ with $M \to \infty$. Here, a recent result of Dolgopyat and Sarig \cite{quenched} shows that for $f(x) = \{x\} - \frac{1}{2}$ and any $x_0 \in \T$, $S_N(f,\alpha,x_0)$ converges (after appropriate centering and scaling) in distribution to a Cauchy random variable. This resembles the behaviour found by Kesten \cite{Kesten1960} who showed that also the \textit{spatial} average 
(that is, $(\alpha,x_0)$ is drawn uniformly at random whereas $N$ is fixed)
converges to a Cauchy distribution.

In the present article, we are dealing with the \textit{quenched temporal} case. This means we are investigating the pointwise behaviour of $S_N(f,\alpha,x_0)$ for fixed $\alpha \in \T$ where we study TDLTs in the sense of \eqref{TLT}.
There are two prominent limit distributions such that a TDLT is satisfied: On the one hand, there are examples where a temporal \textit{central} limit theorem (TCLT) holds, that is, \eqref{TLT} is obtained with $Y$ being a standard Gaussian random variable.
Such results are known to hold for irrational circle rotations for specific irrationals $\alpha$, starting points $x_0 $ and certain functions $f$. For quadratic irrationals $\alpha$, the existence of a TCLT was shown to hold 
for $S_N(f,\alpha,0)$ when $f(x) = \{x\} - 1/2$, $f(x) = \mathds{1}_{[0,\beta)}(x) - \beta$, $\beta \in \mathbb{Q}$ or $f(x) = \log\left\lvert 2 \sin (\pi x)\right\rvert$ (see \cite{beck1,beck2,beck3,borda}).
For the special case where $\alpha = [0;a,a,a,\ldots], a \in \N,$ Borda \cite{borda} showed that a TCLT for $S_N(f,\alpha,0)$ holds for any function $f$ of bounded variation.
The case $f(x) = \mathds{1}_{[0,\beta)}(x) - \beta$ was generalized to arbitrary orbits $S_N(f,\alpha,x_0), x_0 \in \mathbb{R}$ by Dolgopyat and Sarig \cite{dol_sar}
and further by Bromberg and Ulcigrai \cite{brom_ulc} to badly approximable $\alpha$ under some Diophantine assumption (with respect to $\alpha$) on  $\beta$.\\

Note that the results on quadratic irrationals mentioned above do not say anything about typical $\alpha$ since the set of badly approximable numbers (and thus in particular, of quadratic irrationals) is a set of Lebesgue measure $0$. So a natural question is whether a TDLT can hold for almost all $\alpha \in \T$ or at least for $\alpha $ in a set of positive measure.

If $f$ is a smooth function, the existence of a TDLT in the metric sense (i.e. for almost all $\alpha \in \T$) is immediately ruled out: If the Fourier coefficients of $f \sim \sum_{n \in \mathbb{Z}} c_n e(nx)$ decay at rate $c_n = O(1/n^2)$
(which holds in particular for $f \in C^2$), then for almost all $\alpha \in \T$ and all $x_0 \in \mathbb{R}$, $S_N(f,\alpha,x_0)$ is bounded (see \cite{quenched,herman}). Therefore, a TDLT cannot hold because the scaling sequence $(B_M)_{M \in \N}$ needs to be unbounded.
Thus, the interesting functions to consider are those that lack smoothness such as functions that have discontinuities or singularities.
Concerning functions with singularity, Borda \cite{borda} ruled out a central limit theorem for 
$S_N(f,\alpha,0)$ for almost every $\alpha$ where $f(x) =\log (\lvert 2 \sin(\pi x) \rvert)$. In this article, however, we are not considering functions with singularities, but piecewise smooth functions with finitely many discontinuities (compare to, e.g., \cite{dol_sar_no,quenched,upper_dens}).




\begin{df}[Piecewise smooth functions]\label{def_fct1}
    We call a function $f: \mathbb{T} \to \mathbb{R}$ with $\int_{\mathbb{T}} f(x)\,\mathrm{d} \mu(x) = 0$ a
    \textit{piecewise smooth function} if there exist $ \nu \geq 1$ and $ \{\gamma_1, \ldots, \gamma_\nu \} \subseteq \T$ with $ 0 \leq \iota(\gamma_1) < \ldots < \iota(\gamma_{ \nu}) < 1$ 
    ($\iota$ denotes the canonical embedding $\T \hookrightarrow [0,1)$, see Section \ref{Sec_preReq})
    such that the following properties hold:
    \begin{itemize}
        \item $f$ is differentiable on $\mathbb{T}\setminus \{\gamma_1,\ldots,\gamma_{\nu}\}$.
        \item $f'$ extends to a function of bounded variation on $\mathbb{T}$.
         \item There exists an $i \in \{1,\ldots,\nu\}$ such that $\lim_{ \delta \rightarrow 0} \left[ f(\gamma_i - \delta) - f(\gamma_i+ \delta) \right] \neq 0$.
    \end{itemize}
\end{df}

In \cite{upper_dens}, the authors examined the maximal oscillation of $S_N(f,\alpha,x_0)$ for $f$ as in Definition \ref{def_fct1} where 
an unexpected sensitivity on the interplay between the number-theoretic properties of $x_0,\gamma_1, \ldots, \gamma_\nu$ and analytic properties of $f$ was discovered.\\
Note that the class of functions from Definition \ref{def_fct1} contains most of the examples mentioned above, such as $f(x) = \{x\} - 1/2$ or $f(x) = \mathds{1}_{[\beta,\gamma]}, \beta,\gamma \in \T$. 
 Returning to the (non)-existance of TCLTs, the best currently known result for general piecewise smooth $f$ was established in \cite{dol_sar_no}:\\

\noindent {\bf Theorem A. } (Dolgopyat, Sarig, 2018). {\it Let $f$ be a piecewise smooth function as in Definition \ref{def_fct1}. Then there exists a set $\mathcal{E} \subset \T \times \T$ of full two-dimensional (Haar) measure such that for all $(\alpha,x_0) \in \mathcal{E}$, $S_N(f,\alpha,x_0)$ does not satisfy a TDLT.} \\

The aim of the present article is to show that the two-dimensional metric setup above is not necessary and a TDLT fails for almost every $\alpha$ and \textit{any} initial point $x_0 \in \T$:

\begin{thm}\label{thm_no_LT}
    Let $f$ be a piecewise smooth function (see Definition \ref{def_fct1}). Then for (Haar-) almost all $ \alpha \in \mathbb{T}$ and for any $x_0 \in \T$ the following holds:
    Let $N$ be uniformly distributed on $\{1,\ldots,M\}$. Then the sequence of random variables
    $\left(\frac{S_{N}(f,T_{\alpha},x_0) - A_M}{B_M}\right)_{M \in \mathbb{N}}$ does not satisfy a distributional limit theorem in the sense of \eqref{TLT}, regardless of how $(B_M)_{M \in \N}$ and $(A_M)_{M \in \N}$ are chosen.
\end{thm}


\begin{rem*}
Theorem \ref{thm_no_LT} reveals that the set $\mathcal{E}$ from Theorem A can be chosen as
$\mathcal{E} = \mathcal{A} \times \mathbb{T}$ where $\mathcal{A}$ has full ($1$-dimensional) Haar measure. The techniques used in the proof of Theorem A in \cite{dol_sar_no} only allow to make a statement about almost all pairs $(\alpha,x_0) \in \T \times \T$ and we do not know whether adapting the method from \cite{dol_sar_no} would allow to rule out the temporal limit theorem for \textit{every} $x_0 \in \T$
and $\alpha$ in a set $\mathcal{A}$ (that does not depend on $x_0$) of full measure.
Our method of proof takes a different approach and we do not use Fourier-analytic methods as it was done in \cite{dol_sar_no,quenched}.
\end{rem*}

For the special case of the sawtooth function $s(x) = \{x\}-\frac{1}{2}$, Dolgopyat and Sarig showed in \cite[Corollary 2.3]{quenched} that for all starting points $x_0 \in \T$, there exists a set $\mathcal{A}_{x_0} \subseteq \T $ with full Haar measure such that for all $\alpha \in \mathcal{A}_{x_0}$, $S_N(s,\alpha, x_0)$ does not satisfy a TDLT. Again, Theorem \ref{thm_no_LT} implies the stronger result that there exists a set $\mathcal{A} \subseteq \T $ of full Haar-measure such that, for all starting points $x_0  \in \T$ and all $\alpha \in \mathcal{A}$, the associated Birkhoff sum $S_N(s, \alpha, x_0)$ does not satisfy a TLDT.

In \cite[Corollary 2.3]{quenched}, Dolgopyat and Sarig were able to identify a certain family of distributions where each member is realized as
a temporal limit along a suitably normalized subsequence of $S_N(s,T_{\alpha},x_0)$.
In the same paper, the authors ask for a better understanding for general functions in the form of Definition \ref{def_fct1}. A comparable family of distributions appears in our method of proof (see \eqref{def_of_g} in Lemma \ref{indicator_lem}) for all functions $f$ in the form of Definition \ref{def_fct1}. For the special case $f = \mathds{1}_{[0,a]}$, Dolgopyat and Sarig 
 \cite{dol_sar} showed that if $N$ is not sampled uniformly from $\{1,\ldots,M\}$, but $N \sim \text{Log}(\{1,\ldots,M\})$, $S_N(\mathds{1}_{[0,a]},\alpha,0)$ does not satisfy a TDLT.
 However, even for the special case $f = \mathds{1}_{[0,a]}$, the result of Theorem \ref{thm_no_LT} was not yet established.\\

 The rest of this paper is organized as follows. In Section \ref{Sec_preReq}, we fix notation and we state all necessary standard results needed to prove Theorem \ref{thm_no_LT}. In Section \ref{Sec_prepLemmas}, we decompose $f$ into a linear combination of the sawtooth function and certain indicator functions (Proposition \ref{LemRepresentationf}). Further, by using the metric theory of continued fractions, we obtain the almost sure existence of infinitely many (unusually) large partial quotients whose corresponding convergent denominator also satisfies additional properties (see Lemma \ref{Lem_seq_kj} and Remark \ref{Rem_evenoddIndices}). A fact that might be of theoretical interest on its own. In Section \ref{SubSec_proofThm}, Lemma \ref{indicator_lem} establishes limit distributions of $S_N(f, \alpha, x_0)$ along certain subsequences of integers. Finally, we conclude the proof of Theorem \ref{thm_no_LT} by showing that there are at least two such limit distributions that do not coincide.

\section{Prerequisites}
\label{Sec_preReq}

\subsection*{Notation}
Given two functions $f,g:(0,\infty)\to \mathbb{R},$ we write $f(t) = O(g(t)), f \ll g$ or $g \gg f$ if
$\limsup_{t\to\infty} \frac{|f(t)|}{|g(t)|} < \infty$. Any  dependence of the value of the limes superior above on potential parameters is denoted by appropriate subscripts. For two sequences $ (a_k)_{k \in \N}$ and $(b_k)_{k \in \N}$ with $b_k \neq 0$ for all $k \in \N$, we write $a_k \sim b_k, k \to \infty$, if $ \lim_{k \to \infty} \frac{a_k}{b_k}=1$.
We denote the characteristic function of a set $A$ by $\mathds{1}_A$ and understand
the value of empty sums as $0$. For $A \subseteq \N$, we define the lower density of $A$ as
$\liminf\limits_{N \to \infty}\frac{1}{N} \# (A \cap [\![1,N]\!])$.\\

To avoid confusion between elements on $\T \simeq {\Large{\sfrac{\R}{\mathbb{Z}}}}$ and on $\R$, we use the following notation:
We write $\iota: \T \hookrightarrow [0,1)$ for the canonical embedding
$x + \mathbb{Z} \mapsto \{x\} := x - \lfloor x \rfloor$
and let
$\lVert x \rVert
:= \min\{\iota(x),1 - \iota(x)\}$ denote the canonical norm on $\T$. We will denote the normalized Haar measure on $\T$ by $\mu$. For $a,b,x \in \T$, we understand $ \mathbb{1}_{[a,b]}(x) $ as $ \mathbb{1}_{[\iota(a),\iota(b)]}(\iota(x))$. For $a \in \T$ and $n \in \N$, we define as usual $na :=\sum_{i=1}^n a$. If $x \in \R$ and $a \in \T$, we understand $x + a $ as $\iota^{-1}(x) + a \in \T $.\\

Let $X,Y$ be two real-valued random variables defined on a common probability space. If $X$ and $Y$ have the same distribution, we write $X \stackrel{d}{=} Y$. If $X$ has the distribution $\mu$ we write $X \sim \mu$. For $a,b \in \R$ with $a<b$, we denote the uniform distribution on $[a,b]$ as $U([a,b])$. When $a,b \in \N_0$ with $a < b$, $U([\![a,b]\!])$ is the (discrete) uniform distribution on $[a,b] \cap \N_0$.

\subsection*{Continued fractions and Koksma's inequality}\label{cf_prop}
In this subsection, we recall several well-known results from the theory of continued fractions which are heavily used in the proof of Theorem \ref{thm_no_LT}. 
For a more detailed background, we refer the reader to classical literature such as \cite{all_shall,rock_sz}.
Every irrational $\alpha \in [0,1)$ has a unique infinite continued fraction expansion denoted by $[0;a_1,a_2,\dots ]$ with convergents $p_k/q_k := [0;a_1, \dots, a_k]$ that satisfy the recursions
\begin{equation*}
 p_{k+1}  = p_{k+1}(\alpha) = a_{k+1}(\alpha) p_k + p_{k-1}, \qquad q_{k+1}  = q_{k+1}(\alpha) = a_{k+1}(\alpha) q_k + q_{k-1}, \quad k \in \N,
\end{equation*}
with initial values $p_0 = 0,\; p_1 = 1,\; q_0 = 1,\; q_1 = a_1$.
For the sake of brevity, we just write $a_k, p_k, q_k$, although these quantities depend on $\alpha$. Note that the convergents $p_k/q_k$ satisfy the inequalities
\begin{equation}
\label{Eq_parity_deltak}
 \frac{1}{(a_{k+1} +2)q_{k}} \leq \delta_k := (-1)^k (q_k \alpha - p_k) \leq 
 \frac{1}{a_{k+1}q_k}, \quad k \geq 1.
\end{equation}






Conversely, if $\lvert \alpha - p/q\rvert < \frac{1}{2q^2}$, Legendre's Theorem implies that $p/q$ is a convergent of $\alpha$.\\

Since this article deals with almost sure behaviour, we also make use of the following classical results that arise from the well-studied area of the \textit{metric} theory of continued fractions:

\begin{itemize}
\item (Diamond and Vaaler \cite{diamond_vaaler}): For almost every $\alpha$,
\begin{equation}
\label{Eq_trimmed_sum}
\sum_{\ell \leq K} a_{\ell} - \max\limits_{\ell \leq K} a_{\ell} \sim \frac{K \log K}{\log 2}, \quad K \to \infty.
\end{equation} 
\item (Khintchine and L\'{e}vy, see, e.g., \cite[Chapter 5, §9, Theorem 1]{rock_sz}): For almost every $\alpha$,
\begin{equation}
\label{Eq_size_of_q_k}
\log q_k \sim \tfrac{\pi^2}{12 \log 2} k, \quad  k \to \infty.
\end{equation}
\end{itemize}


On several positions in the proof, we will make use of Koksma's inequality which allows to estimate the error between sums and corresponding integrals. For more details about this topic and the closely related area of Discrepancy theory, we refer the reader to \cite{kuipers}.
Denoting the discrepancy of a sequence $\yn \subseteq \T$ at stage $N \in \N$ by
\[D_N(\yn) := \sup\limits_{0 \leq a \leq b < 1} \left \lvert \frac{1}{N}\#\left\{1 \leq n \leq N: \iota(y_n) \in [a,b]\right\} - (b-a) \right \rvert \]
and the total variation of $f: \T \to \R$ by $\Var(f)$,
Koksma's inequality is given by

\[\left\lvert \sum_{i=1}^{N} f(y_i) - N\int_{\T} f(x) \mathrm{d}\mu(x) \right\rvert \leq \Var(f)N D_N(\yn).\]

In the special case where $\yn$ is the Kronecker sequence $(n\alpha)_{n \in \N}$, we have the estimates

\[D_{q_n}(\yn) \ll \frac{1}{q_n}, \quad D_{N}(\yn) \ll \frac{1}{N}\sum_{i=1}^{k} a_i,\]

where $k = k(N)$ is such that $q_{k-1} \leq N < q_{k}$. Thus Koksma's inequality leads (in this particular case also known as Denjoy-Koksma inequality, see, e.g., \cite{herman}) to 

\begin{equation}\label{denj_koks} \left \lvert S_{q_n}(f, \alpha, x_0) \right \rvert
\ll_f 1, \quad
\left \lvert S_{N}(f, \alpha, x_0) \right \rvert \ll_f \sum_{i=1}^{k} a_i,
\end{equation}

with the implied constant being uniform in $x_0$.

\section{Proof of Theorem \ref{thm_no_LT}}
\label{Sect_Pf_noCLT}
\subsection{Preparatory Lemmas}
\label{Sec_prepLemmas}

\begin{proposition}\label{LemRepresentationf}
Let $f : \T \rightarrow \R$ be as in Definition \ref{def_fct1}.

Let $h:\T \to \R$ be defined as \[h(x) = \sum_{i=1}^{\nu} H_i\left ( \iota(x) - \frac{1}{2} \right)  + \sum_{i=1}^{\nu} H_i \left(  \mathbb{1}_{[0, \gamma_i)}(x) -\iota(\gamma_i) \right),\]
where 
$H_i := \lim_{ \delta \rightarrow 0} \left[f(\gamma_i - \delta) - f(\gamma_i + \delta) \right]$.
Then, for almost every $\alpha \in \T$, any $N \in \N$ and any $y \in \T$, we have
\[S_N(f,\alpha,y) = S_N(h,\alpha,y) + O_f(1),\]
with the implied constant only depending on $f$.
\end{proposition}

\begin{proof}
    This can be proven analogously to \cite[Lemma 3.1]{upper_dens}. A more detailed proof can be found in \cite[Appendix A]{quenched}.\end{proof}

\begin{proposition}\label{ds_thm} (Duffin and Schaeffer, \cite[Theorem 3]{duff_sch}).
Let $A \subseteq \N$  be a set of positive lower density and $\psi: \N \to [0,\infty)$ be a monotone decreasing function such that $\sum\limits_{q =1}^{\infty} \psi(q) = \infty$. Then, for almost every $\alpha$, there exist infinitely many coprime $(p,q) \in \mathbb{Z} \times A$ that satisfy $\left\lvert \alpha - \frac{p}{q}\right\rvert < \frac{\psi(q)}{q}$.
\end{proposition}

\begin{proposition}\label{gallagher}(Gallagher, \cite[Lemma 2]{gall}).
Let $(I_k)_{k \in \N} \subseteq \T$ be a sequence of intervals with $\lim_{k \to \infty} \mu(I_k) = 0$.
Further let $c > 0$ and $(U_k)_{k \in \N}$ be a sequence of measurable sets that satisfy the following for all $k \in \N$:  
\begin{itemize}
    \item $U_k \subseteq I_k$,
    \item $\mu(U_k) \geq c \mu(I_k)$.
\end{itemize}
Then, $\mu(\limsup\limits_{k \to \infty} U_k) = \mu(\limsup\limits_{k \to \infty} I_k)$.
\end{proposition}

Combining the statements above, we can deduce the following result.

\begin{lem}
\label{Lem_seq_kj}
Let $A \subseteq \N$ be a set with positive lower density. Then, for almost every $\alpha = [\mathbb{Z};a_1, a_2, \ldots] \in \T$, there exists a sequence of even integers $(k_j)_{j \in \N}$ such that $q_{k_j} \in A$ for all $j \in \N$ and $ \lim_{j \to \infty} \frac{\sum_{i=1}^{k_j} a_i }{a_{k_j +1}} = 0$.
\end{lem}
\begin{proof}
Let $ \psi(q) = \frac{1}{q \log q \log \log q \log \log \log q }$ \footnote{For convenience, we set $ \log x := 1$ if $x \leq e$.}, then it holds that $ \sum_{q \in \N } \psi(q)= \infty$ as well as $\psi(q) \leq 1$ for all $q \in \N$.
\par{}
Let $(r_k/s_k)_{k \in \N}$ be the set of rationals with $s_k \in A$ and $1 \leq r_k \leq s_{k}-1$ with $\gcd(r_k,s_k) =1$. We define 
\begin{equation*}
I_k:= \iota^{-1}\left( \frac{r_k}{s_k} - \frac{\psi(s_k)}{s_k}, \frac{r_k}{s_k} + \frac{\psi(s_k)}{s_k} \right) \quad \text{and}  \quad U_k := \iota^{-1}\left( 0, \frac{r_k}{s_k} + \frac{\psi(s_k)}{s_k} \right).
\end{equation*}
By Proposition \ref{ds_thm}, 
we have $\mu(\limsup_{k \to \infty} I_k) = 1$.
 Since clearly $U_k \subseteq I_k$ and $ \mu( U_k) \geq \frac{1}{2}\mu(I_k)$ for all $k \in \N$, an application of Proposition \ref{gallagher} shows 
 $ \mu( \limsup_{k \rightarrow \infty} U_k) = 1$. In other words, for almost all $\alpha \in \T$, 
there are infinitely many coprime pairs $(p,q) \in \N \times A$ such that 
\begin{equation}
\label{EqApproxpq}
    0 \leq \alpha - \frac{p}{q} < \frac{\psi(q)}{q} = \frac{1}{q^2 \log q \log \log q \log \log \log q}. 
\end{equation}
By Legendre's Theorem, for $q \geq 10$, the above is only possible if $ (p,q)$ is a convergent of $\alpha $. Thus, the pairs $(p,q), q \geq 10$ that satisfy \eqref{EqApproxpq} form a subsequence $ (p_{k_j}, q_{k_j})_{j \in \N}$ of the sequence of convergents $ ( p_k, q_k)_{k \in \N}$. Since $\alpha - \frac{p_{k_j}}{q_{k_j}} \geq 0 $ for all $j \in \N$, it follows by \eqref{Eq_parity_deltak} that all $k_j$ are even. Moreover, by construction of $\psi$ and combining  \eqref{Eq_parity_deltak} and \eqref{Eq_size_of_q_k},
we have $ a_{k_j+1} \gg k_j \log k_j \log \log k_j$. By \eqref{Eq_trimmed_sum} this implies that for almost every $\alpha$, we have $\sum\limits_{i=1}^{k_j} a_i =  o \left( a_{k_j +1} \right)   $.
\end{proof}

\begin{rem}
\label{Rem_evenoddIndices}
By obvious modifications, the statement of Lemma \ref{Lem_seq_kj} also holds when ``even'' is replaced by ``odd''.
In Lemma \ref{indicator_lem}, this would lead to an even larger class of limiting distributions that are realized as limits of certain Birkhoff sums along suitable subsequences.
For our purpose of ruling out any TDLT, the stated version of Lemma \ref{Lem_seq_kj} is sufficient.
\end{rem}

\begin{proposition}\label{isolated_density}
Let $\beta_1,\beta_2,\ldots,\beta_{\nu} \in \T \setminus \{0\}, \nu \in \N$. Then there exists $\delta> 0$ such that the set $\{N \in \N: \forall 1 \leq j \leq \nu: \lVert N\beta_j \rVert > \delta \}$ has positive lower density.
\end{proposition}

\begin{proof}
We partition $\{ \beta_i \}_{i=1}^{\nu} $ into rational and irrational numbers. Without loss of generality, we may assume $\iota(\beta_1) = \frac{a_1}{b_1},\ldots,\iota(\beta_k) = \frac{a_k}{b_k} \in \Q$ with $a_i,b_i \in \N, \gcd(a_i,b_i) = 1, b_i \geq 2$ since $\beta_i \neq 0$ for $i = 1,\ldots,k$, and 
$\iota(\beta_{k+1}),\ldots,\iota(\beta_{\nu}) \notin \Q$. Let $b_{\pi} := \prod_{i=1}^{k}b_i$.
Clearly, if $N \equiv 1 \pmod {b_{\pi}}$, then for all $1 \leq i \leq k$, $b_i \nmid N$ and thus,
$\iota(N\beta_i) \in \left\{\frac{1}{b_i},\ldots,\frac{b_i-1}{b_i}\right\}$, which is disjoint from 
$(0,\delta) \cup (1-\delta,1)$ if $\delta$ is chosen sufficiently small. Since $\{N \in \N: N \equiv 1 \pmod{b_{\pi}}\}$ has positive lower density, it suffices to show that
\begin{equation*}
 \left\{M \in \mathbb{N}:   \forall i \in \{k+1,\ldots,\nu\}: 
 \lVert (M b_{\pi}+1)\beta_i \rVert > \delta \right\}
\end{equation*}

has positive lower density. Since $ \iota ( b_{\pi} \beta_i ) \notin \Q$ for all $i=k+1, \ldots, \nu$, it follows that $ \left \{ (M b_{\pi} \beta_i + \beta_i )\right\}_{M \in \N}$ is uniformly distributed on $\T$. This immediately shows

\[\liminf_{N \to \infty}
\frac{1}{N}\#\left\{M \leq N:  \forall i \in \{k+1,\ldots,\nu\}: 
 \lVert (M b_{\pi}+1)\beta_i \rVert  > \delta \right\}
 \geq 
 1 - 2\nu\delta > 0,
\]
 provided $\delta < \frac{1}{2\nu}$.    
\end{proof}

\subsection{Main Lemma and conclusion of the proof}
\label{SubSec_proofThm}
\begin{lem}\label{indicator_lem}
Let $f(x) =\left(\sum_{i=1}^{\nu}  H_i\right) \left(\iota(x) - \frac{1}{2}\right) +  \sum_{i=1}^{\nu}  H_i \left( \mathbb{1}_{[0, \gamma_i)} (x) - \iota(\gamma_i) \right)$
where $\gamma_1,\ldots,\gamma_{\nu} \in \mathbb{T}$ are distinct. Then for almost every $\alpha = [\mathbb{Z};a_1,a_2,\ldots] \in \T$ and any $x_0 \in \T$, there exists 
an increasing sequence $({\kj})_{\ell \in \N}$ such that the following holds:

    \begin{itemize}
    \item For every $\ell \in \N$, $q_{\kj}$ is a denominator of a convergent of $\alpha$.
    \item $\lim\limits_{\ell \to \infty} \frac{\sum\limits_{i =1}^{\kj} a_i }{a_{\kj+1}} = 0$.
    \item The limits \[\overline{x_0} := \lim_{\ell \to \infty} q_{\kj}x_0, \quad \overline{\gamma_i} := \lim_{\ell \to \infty} q_{\kj}\gamma_i, \quad i = 1,\ldots,\nu \quad\] exist and satisfy
    $\overline{\gamma_i} \neq \overline{\gamma_1}$ for all $i = 2,\ldots, \nu$.
    If $\gamma_1 \neq 0$, then $\overline{\gamma_1} \neq 0$.
    \item We have
    \begin{equation}
    \label{def_of_g}
    \lim_{\ell \to \infty}
    \sup_{c \in [0,1]} \left\lvert
    \frac{S_{ \lfloor c a_{\kj+1}\rfloor q_{\kj} }(f,\alpha,x_0)}{a_{\kj+1}} -  \left( \left( \sum_{i = 1}^{\nu} H_i \right) \left(\int_{0}^{c} \iota(y + \overline{x_0}) \, \mathrm{d}y - \frac{c}{2}\right) + \sum_{i = 1}^{\nu} H_i\left(\int_{0}^{c} \mathds{1}_{\left[0, \overline{\gamma_i}\right]}\left(y + \overline{x_0}\right) \,\mathrm{d}y - \iota(c\overline{\gamma_i})\right)\right)\right\rvert
    = 0.
    \end{equation}
    \end{itemize}
\end{lem}

\begin{proof}
Without loss of generality, we may assume that $\gamma_i \neq 0$ for all $i = 1,\ldots\nu$, since otherwise, we observe that $S_N(f,\alpha,x_0) = S_N(\tilde{f},\alpha,x_0 + y_0)$ where $\tilde{f}(x) = f(x- y_0)$ and $y_0$ is chosen such that $\gamma_i + y_0 \neq 0$ for all $i = 1,\ldots,\nu$. 
We apply Proposition \ref{isolated_density} to
\[\left\{\beta_1,\ldots,\beta_{2\nu-1}\right\} := \left\{\gamma_1,\ldots,\gamma_{\nu},\gamma_2 - \gamma_1,\gamma_3-\gamma_1,\ldots\gamma_{\nu}-\gamma_{1}\right\},\]


which gives us a set $A \subseteq \N$ with positive lower density and $\delta > 0$ such that for all $N \in A$,
$ \lVert N\gamma_i\rVert > \delta$ and $\lVert N\gamma_i - N\gamma_1\rVert > \delta$ for $ i = 2,\ldots,\nu$.

Next, we apply Lemma \ref{Lem_seq_kj} to $A$:
For almost every $\alpha \in \T$, there exists a sequence $(q_{k_j})_{j \in \N} \subseteq A$ such that the following holds:

   \begin{itemize}
    \item $k_j$ is even for all $j \in \N$.
    \item For every $j \in \N$, $q_{k_j}$ is a denominator of a convergent of $\alpha$.
    \item $\lim\limits_{j \to \infty} \frac{\sum\limits_{i =1}^{k_j} a_i }{a_{k_j+1}} = 0$.
    \end{itemize}

For fixed $x_0 \in \T$, observe that $(q_{k_j}x_0)_{j \in \N},\; (q_{k_j}\gamma_i)_{j \in \N}, \; i = 1,\ldots,\nu$ are (bounded) sequences in $\T$, thus there exists a subsequence $(n_{\ell})_{\ell \in \N}$ of $(k_{j})_{j \in \N} $ such that the limits
\[\overline{x_0} := \lim_{\ell \to \infty} q_{\kj}x_0, \quad \overline{\gamma_i} := \lim_{\ell \to \infty} q_{\kj}\gamma_i, \quad i = 1,\ldots,\nu \quad\] all exist. Since $(q_{\kj})_{\ell \in \N} \subseteq (q_{k_j})_{j \in \N} \subseteq A$, we have for any $\ell \in \N$, $\lVert q_{\kj}\gamma_i - q_{\kj}\gamma_1\rVert > \delta$ and thus, 
$\overline{\gamma_i} \neq \overline{\gamma_1}$ for all $i = 2,\ldots,\nu$.

We now turn our attention to prove \eqref{def_of_g}.
First, we prove that for $s(x) = \iota(x) - \frac{1}{2}$, we have
\begin{equation}
\label{sawtooth_fct}
\lim_{\ell \to \infty} \frac{S_{ \lfloor c a_{\kj+1}\rfloor q_{\kj}}(s,\alpha,x_0)}{a_{\kj+1}} =  \int_{0}^{c} \iota(y + \overline{x_0}) \, \mathrm{d}y - \frac{c}{2},\end{equation}
with the convergence being uniform in $c \in [0,1]$.
Let $\varepsilon > 0$ be given. 
We will show that for any sufficiently large $\ell$ and any integer $u$ with
$0 \leq u \leq \lfloor c a_{\kj+1}\rfloor$ that satisfies
$\left\lVert\frac{u}{a_{\kj +1}} + \overline{x_0}\right\rVert > \varepsilon$,
we have

\begin{equation}\label{partial_sum_sawtooth}\left\lvert\left(S_{(u+1)q_{\kj}}(s,\alpha,x_0) - S_{uq_{\kj}}(s,\alpha,x_0)\right)
 -  \left \{ \frac{u}{a_{\kj +1}} + \iota(\overline{x_0}) \right \} - \frac{1}{2} \right\rvert < \varepsilon.
\end{equation}

For $\ell$ large enough, we have
\[    \left\lVert q_{\kj} x_0  - \overline{x_0}\right\rVert < \varepsilon/10.
\]

Now observe that   
\[
\begin{split}
S_{(u+1)q_{\kj}}(s,\alpha,x_0) - S_{uq_{\kj}}(s,\alpha,x_0)
&= S_{q_{\kj}}\left(s, \alpha, T^{uq_{\kj}}(x_0)\right) \\
&=\sum_{n = 0}^{q_{\kj}-1}  \left \{ \iota\left ( (n + uq_{\kj})\alpha \right) + \iota(x_0)\right \} - \frac{q_{\kj}}{2}\\
&= 
\sum_{n = 0}^{q_{\kj}-1}  \left \{ n\frac{p_{\kj}}{q_{\kj}} + n\frac{\delta_{\kj}}{q_{\kj}}+u\delta_{\kj} + \iota(x_0) \right\}
- \frac{q_{\kj}}{2} \\
&= \sum_{n = 0}^{q_{\kj}-1}  \left\{ n\frac{p_{\kj}}{q_{\kj}} +\frac{u/a_{\kj+1}}{q_{\kj}} + \frac{O(1/a_{\kj+1})}{q_{\kj}} + \iota(x_0) \right\} 
- \frac{q_{\kj}}{2}
\end{split}
\]
where $\delta_{\kj} := \iota( q_{\kj}  \alpha  )  = \frac{1}{a_{\kj+1}q_{\kj}}\left(1 + O\left(\frac{1}{a_{\kj+1}}\right)\right)$, which follows from the assumption that $\kj$ is even and we apply  \eqref{Eq_parity_deltak}. Since $\gcd(p_{\kj},q_{\kj}) =1$, we have 

\[\begin{split}
&\sum_{n = 0}^{q_{\kj}-1}   \left \{ n\frac{p_{\kj}}{q_{\kj}} +\frac{u/a_{\kj+1}}{q_{\kj}} + \frac{O(1/a_{\kj+1})}{q_{\kj}} + \iota(x_0) \right\} \\
=&
\sum_{j = 0}^{q_{\kj}-1} \left \{ \frac{j}{q_{\kj}} +\frac{u/a_{\kj+1}}{q_{\kj}} + \frac{O(1/a_{\kj+1})}{q_{\kj}} + \frac{\lfloor q_{\kj} \iota( x_0)\rfloor}{q_{\kj}} + \frac{\iota(q_{\kj} x_0)}{q_{\kj}} \right \}\\
=&
\sum_{j = 0}^{q_{\kj}-1} \left\{ \frac{j}{q_{\kj}} +\frac{u/a_{\kj+1} + \iota \left(q_{\kj}  x_0\right)}{q_{\kj}} + \frac{O(1/a_{\kj+1})}{q_{\kj}}\right\}
\\
=&
\sum_{j = 0}^{q_{\kj}-1}  \left\{ \frac{j}{q_{\kj}} +\frac{u/a_{\kj+1} + \iota \left(\overline{x_0}\right)}{q_{\kj}} + \frac{O(1/a_{\kj+1})}{q_{\kj}} + \frac{R_\varepsilon}{q_{\kj}}\right\},
\end{split}
\]
where $ R_{\varepsilon} := \iota( q_{\kj}x_0) - \iota( \overline{x_0})$ which satisfies $\lvert R_{\varepsilon} \rvert \leq \frac{\varepsilon}{10}$ by the choice of $\ell$. For all integers $u$ with $0 \leq u \leq \lfloor c a_{\kj +1} \rfloor $ such that $ \left\lVert \frac{u}{a_{\kj +1}} +  \overline{x_0}  \right\rVert > \varepsilon$, we have
\begin{align*}
S_{(u+1)q_{\kj}}(s,\alpha,x_0) - S_{uq_{\kj}}(s,\alpha,x_0)
&= S_{q_{\kj}}\left(s, \alpha, T_{\alpha}^{u q_{\kj}}(x_0)\right) \\ 
&= \sum_{j = 0}^{q_{\kj}-1}  \left\{ \frac{j}{q_{\kj}} +\frac{u/a_{\kj+1} + \iota \left(\overline{x_0}\right)}{q_{\kj}} + \frac{O(1/a_{\kj+1})}{q_{\kj}} + \frac{R_\varepsilon}{q_{\kj}}\right\} - \frac{q_{\kj}}{2} \\
&= \sum_{j = 0}^{q_{\kj}-1} \left(  \frac{j}{q_{\kj}} +\frac{ \left \{u/a_{\kj+1} + \iota \left(\overline{x_0}\right) \right \}}{q_{\kj}} + \frac{O(1/a_{\kj+1})}{q_{\kj}} + \frac{R_\varepsilon}{q_{\kj}} \right) - \frac{q_{\kj}}{2}\\
& = \left\{u/a_{\kj+1} + \iota(\overline{x_0}) \right\} - \frac{1}{2} + O(1/a_{\kj+1}) + R_{\varepsilon},
\end{align*}
which proves \eqref{partial_sum_sawtooth}. Clearly, 

\[\#\left\{
0 \leq u \leq \lfloor c a_{\kj+1}\rfloor:
\left\lVert\frac{u}{a_{\kj +1}} + \overline{x_0}\right\rVert < \varepsilon
\right\} \leq 2\varepsilon a_{\kj+1} + 2\]

and by the Denjoy-Koksma inequality (see \eqref{denj_koks}), we have

\[\lvert S_{(u+1)q_{\kj}}(s,\alpha,x_0) - S_{uq_{\kj}}(s,\alpha,x_0)\rvert
\ll 1,
\]
for any $0 \leq u \leq {a_{\kj+1}-1}$.
Thus,  
\[\begin{split}S_{ \lfloor c a_{\kj+1}\rfloor}q_{\kj}(s,\alpha,x_0)
&= \sum_{u = 0}^{\lfloor c a_{\kj+1}\rfloor-1} S_{(u+1)q_{\kj}}(s,\alpha,x_0) - S_{uq_{\kj}}(s,\alpha,x_0)
\\&= \sum_{u = 0}^{\lfloor c a_{\kj+1}\rfloor-1}\left(\left\{u/a_{\kj+1} + \iota(\overline{x_0})\right\} - \frac{1}{2} + O(\varepsilon) + O\left(1/a_{\kj+1}\right)\right)
+ O\left(\varepsilon a_{\kj+1}\right)
\\
&= a_{\kj+1}\left(\int_{0}^{c} \iota\left(y + \overline{x_0}\right) \,\mathrm{d}y - \frac{c}{2} + O(\varepsilon) \right) + O(1),
    \end{split}
\]
where the implied constants in the $O$-terms depend neither on $c$ nor on $\varepsilon$.
In the last line, we used Koksma's inequality to compare sum and integral. With $\varepsilon \to 0$, \eqref{sawtooth_fct} follows.

Next, we fix $i \in \{1,\ldots,\nu\}$ and show that
 \begin{equation}
 \label{single_indicator}
 \lim_{\ell \to \infty} \frac{S_{ \lfloor c a_{\kj+1}\rfloor q_{\kj}}(\mathds{1}_{[0,\gamma_i]},\alpha,x_0)}{a_{\kj+1}} =  \left(\int_{0}^{c} \mathds{1}_{\left[0, \overline{\gamma_i}\right]}\left(y + \overline{x_0}\right) \,\mathrm{d}y - c\iota(\overline{\gamma_i})\right),\end{equation}
 with the convergence being uniform in $c \in [0,1]$.
   For convenience, we will drop the index $i$ in the following, that is, we set 
   $\gamma := \gamma_i, \overline{\gamma} := \overline{\gamma_i}$. Since $\overline{\gamma} \neq 0$, we take $\varepsilon > 0$ such that $\varepsilon \leq \lVert \overline{\gamma}\rVert$.
    Further, let $\ell$ be large enough such that
    $\left\lvert \frac{\lfloor c a_{\kj+1}\rfloor}{a_{\kj+1}} - c\right\rvert = \frac{\{ca_{\kj+1}\}}{a_{\kj+1}}< \varepsilon/10$ uniformly in $c \in [0,1]$. Moreover for $\ell$ large enough, we have
    \[\left\lVert q_{\kj} \gamma - \overline{\gamma}\right\rVert < \varepsilon/10, \quad
    \left\lVert q_{\kj} x_0  - \overline{x_0}\right\rVert < \varepsilon/10.
    \]

We will show that for any $\ell$ sufficiently large and any 
$0 \leq u \leq \lfloor c a_{\kj+1}\rfloor$ that satisfies
$\left\lVert\frac{u}{a_{\kj+1}} + \overline{x_0} - \overline{\gamma}\right\rVert > \varepsilon$
and 
$\left\lVert\frac{u}{a_{\kj+1}} + \overline{x_0}\right\rVert > \varepsilon$,
we have

    \[\left\lvert\left(S_{(u+1)q_{\kj}}(\mathds{1}_{[0,\gamma]},\alpha,x_0) - S_{uq_{\kj}}(\mathds{1}_{[0,\gamma]},\alpha,x_0)\right)
    - \left(\mathds{1}_{\left[ 0,\overline{\gamma} \right]} \left( u/a_{\kj+1} + \overline{x_0}  \right)
    -  \iota(\overline{\gamma})\right)\right\rvert < \varepsilon.
    \]

To prove this, observe that

\[\begin{split}
S_{(u+1)q_{\kj}}(\mathds{1}_{[0,\gamma]},\alpha,x_0) - S_{uq_{\kj}}(\mathds{1}_{[0,\gamma]},\alpha,x_0)
&= S_{q_{\kj}}\left(\mathds{1}_{[0,\gamma]}, \alpha, T_{\alpha}^{u q_{\kj}}(x_0)\right) \\ 
&= \#\left\{0 \leq n \leq q_{\kj}-1: \iota(n \alpha +uq_{\kj}\alpha + x_0) \in [0,\iota(\gamma)] \right\}
- \iota(\gamma)q_{\kj}\\
&= \#\left\{0 \leq n \leq q_{\kj}-1: \iota \left( n\frac{p_{\kj}}{q_{\kj}} + n\frac{\delta_{\kj}}{q_{\kj}}+u\delta_{\kj} + x_0 \right) \in [0,\iota(\gamma)]\right\}
- \iota(\gamma) q_{\kj}\\
&= \#\left\{0 \leq n \leq q_{\kj}-1: \iota \left( n\frac{p_{\kj}}{q_{\kj}} +\frac{u/a_{\kj+1}}{q_{\kj}} + \frac{O(1/a_{\kj+1})}{q_{\kj}} + x_0 \right) \in [0,\iota(\gamma)]\right\} \\
& \qquad - \iota(\gamma) q_{\kj}.
\end{split}\]
Since $\gcd(p_{\kj},q_{\kj}) =1$, we have

\[\begin{split}
& \phantom{=}  \#\left\{0 \leq n \leq q_{\kj}-1: \iota \left( n\frac{p_{\kj}}{q_{\kj}} +\frac{u/a_{\kj+1}}{q_{\kj}} + \frac{O(1/a_{\kj+1})}{q_{\kj}} + x_0 \right) \in [0, \iota(\gamma)]\right\}
\\
& = \#\left\{0 \leq n \leq q_{\kj}-1 : \left \{ n \frac{p_{\kj}}{q_{\kj}} +\frac{u/a_{\kj+1}}{q_{\kj}} + \frac{O(1/a_{\kj+1})}{q_{\kj}} + \frac{\lfloor q_{\kj} \iota( x_0) \rfloor}{q_{\kj}} + \frac{\iota(q_{\kj} x_0)}{q_{\kj}} \right \} \in \left[0,\frac{\lfloor q_{\kj} \iota(\gamma) \rfloor}{q_{\kj}} + \frac{ \iota( q_{\kj}\gamma )}{q_{\kj}}\right]\right\} \\
&= \#\left\{0 \leq j \leq q_{\kj}-1: \left \{ \frac{j}{q_{\kj}} +\frac{ \left \{ u/a_{\kj+1} + \iota( \overline{x_0} ) \right \}  }{q_{\kj}} + \frac{O(1/a_{\kj+1})}{q_{\kj}} + \frac{R_{\varepsilon}}{q_{\kj}} \right \} \in \left[0,\frac{\lfloor q_{\kj} \iota(\gamma) \rfloor}{q_{\kj}} + \frac{ \iota( q_{\kj}\gamma )}{q_{\kj}}\right]\right\}\\
&= \#\left\{0 \leq j \leq q_{\kj}-1: \left \{ \frac{j}{q_{\kj}} +\frac{ \left \{ u/a_{\kj+1} + \iota(\overline{x_0}) \right \} }{q_{\kj}} + \frac{O(1/a_{\kj+1})}{q_{\kj}} + \frac{R_\varepsilon}{q_{\kj}} \right \}  \in \left[0,\frac{\lfloor q_{\kj}\iota(\gamma) \rfloor}{q_{\kj}} + 
\frac{\iota(\overline{\gamma})}{q_{\kj}} +
\frac{S_\varepsilon}{q_{\kj}}\right]\right\},
\end{split}
\]

where $ R_{\varepsilon}:= \iota( q_{n_\ell} x_0) - \iota( \overline{x_0}) $ and $S_{\varepsilon} := \iota( q_{n_\ell} \gamma ) - \iota( \overline{\gamma}) $. By the choice of $\ell$, we have $ \lvert  R_{\varepsilon} \rvert, \lvert S_{\varepsilon} \rvert \leq \frac{\varepsilon}{10} \leq \frac{\left\lVert\overline{\gamma}\right\rVert}{10}$. Moreover, let $\ell$ be large enough such that
$\lvert O(1/a_{\kj}))\rvert \leq \frac{\varepsilon}{10}$ and since $\left\lVert\frac{u}{a_{\kj +1}} + \overline{x_0}\right\rVert > \varepsilon$, we can drop the fractional part in the previous expression.

We now distinguish two cases: 
First, consider $u$ with $\left\{u/a_{\kj+1} + \iota(\overline{x_0})\right\} \leq \iota(\overline{\gamma})$. Then using our assumption
$\left\lVert\frac{u}{a_{\kj +1}} + \overline{x_0} - \overline{\gamma}\right\rVert > \varepsilon$, it follows that $ \left\{ u/a_{\kj +1} + \iota(\overline{x_0})  \right\} - \iota( \overline{\gamma} ) \leq - \varepsilon$ and thus

\[0 \leq \left\{u/a_{\kj+1} + \iota(\overline{x_0})\right\} + R_{\varepsilon} + O(1/a_{\kj+1}) - S_{\varepsilon} < \iota( \overline{\gamma}) .\]

This implies

\begin{equation}
\label{case_smaller}
\#\left\{0 \leq j \leq q_{\kj}-1: \frac{j}{q_{\kj}} +\frac{\left\{u/a_{\kj+1} + \iota(\overline{x_0})\right\}}{q_{\kj}} + \frac{O(1/a_{\kj+1})}{q_{\kj}} + \frac{R_\varepsilon}{q_{\kj}} \in \left[0,\frac{\lfloor q_{\kj} \iota(\gamma) \rfloor}{q_{\kj}} + 
\frac{\iota(\overline{\gamma})}{q_{\kj}} +
\frac{S_\varepsilon}{q_{\kj}}\right]\right\} = \lfloor q_{\kj}\iota(\gamma) \rfloor +1.
\end{equation}

Similarly, if $\left\{u/a_{\kj+1} + \iota(\overline{x_0})\right\} > \iota(\overline{\gamma})$, then

\begin{equation}
\label{case_larger}
\#\left\{0 \leq j \leq q_{\kj}-1: \frac{j}{q_{\kj}} +\frac{\left\{u/a_{\kj+1} + \iota(\overline{x_0})\right\}}{q_{\kj}} + \frac{O(1/a_{\kj+1})}{q_{\kj}} + \frac{R_\varepsilon}{q_{\kj}} \in \left[0,\frac{\lfloor q_{\kj}\iota(\gamma) \rfloor}{q_{\kj}} + 
\frac{\iota(\overline{\gamma})}{q_{\kj}} +
\frac{S_\varepsilon}{q_{\kj}}\right]\right\} = \lfloor q_{\kj} \iota(\gamma) \rfloor.
\end{equation}

Note that

\[\iota( \gamma )q_{\kj} = \lfloor q_{\kj} \iota(\gamma) \rfloor + \iota(\overline{\gamma}) + T_{\varepsilon},\]
where $T_{\varepsilon} = \iota( q_{\kj} \gamma) - \iota( \overline{\gamma})$. We choose $\ell$ large enough such that $\lvert T_{\varepsilon}\rvert \leq \frac{\varepsilon}{10}$. Combining this with \eqref{case_smaller} and \eqref{case_larger} yields
\begin{align*}
\left\lvert\left(S_{(u+1)q_{\kj}}(\mathds{1}_{[0,\gamma]},\alpha,x_0) - S_{uq_{\kj}}(\mathds{1}_{[0,\gamma]},\alpha,x_0)\right)
- \left(\mathds{1}_{\left[
 0, \overline{\gamma} \right]} \left(  u/a_{\kj+1} + \overline{x_0} \right)
-  \iota(\overline{\gamma})\right)\right\rvert = \lvert T_{\varepsilon} \rvert< \varepsilon,
\end{align*}
 for any $0 \leq u \leq \lfloor c a_{\kj+1}\rfloor$ that satisfies
$\left\lVert\frac{u}{a_{\kj +1}} + \overline{x_0} - \overline{\gamma}\right\rVert > \varepsilon$
and 
$\left\lVert\frac{u}{a_{\kj+1}} + \overline{x_0}\right\rVert > \varepsilon$. Clearly, 

\[\#\left\{
0 \leq u \leq \lfloor c a_{\kj+1}\rfloor:
\left\lVert\frac{u}{a_{\kj+1}} + \overline{x_0} - \overline{\gamma}\right\rVert \leq \varepsilon
\text{ or } \left\lVert\frac{u}{a_{\kj+1}} + \overline{x_0}\right\rVert \leq \varepsilon
\right\} \leq 4\varepsilon a_{\kj+1} + 4\]


and thus analogously to above, we obtain
\[
\begin{split}
S_{ \lfloor c a_{\kj+1}\rfloor}q_{\kj}(\mathds{1}_{[0,\gamma]},\alpha,x_0)
&= \sum_{u = 0}^{\lfloor c a_{\kj+1}\rfloor-1} S_{(u+1)q_{\kj}}(\mathds{1}_{[0,\gamma]},\alpha,x_0) - S_{uq_{\kj}}(\mathds{1}_{[0,\gamma]},\alpha,x_0) \\
&= \sum_{u = 0}^{\lfloor c a_{\kj+1}\rfloor-1}\left(\mathds{1}_{ \left[ 0,  \overline{\gamma} \right] } \left( 
 u/a_{\kj+1} + \overline{x_0} \right) 
-  \iota(\overline{\gamma})\right) + O(\varepsilon a_{\kj +1}) + O(1)\\
&= a_{\kj+1}\left(\int_{0}^{c} \left( \mathds{1}_{ \left[ 0,  \overline{\gamma} \right] } \left( 
y + \overline{x_0} \right) 
-  \iota(\overline{\gamma}) \right)\mathrm{d}y\right) + O(\varepsilon a_{\kj +1}) + O(1).
\end{split}
\]
With $\varepsilon \to 0$, \eqref{single_indicator} follows. Combining \eqref{sawtooth_fct} and \eqref{single_indicator}, we obtain statement \eqref{def_of_g}, which finishes the proof.
\end{proof}

\begin{proof}[Proof of Theorem \ref{thm_no_LT}]
We assume that there exist normalizing sequences $(A_M)_{M \in \N}$ and $(B_M)_{M \in \N}$ with $A_M \in \R, \; B_M > 0$ and $ B_M \to \infty$ such that
\begin{equation}
\label{LiminDist}
\lim_{M \rightarrow \infty} \frac{S_N(f, \alpha, x_0) - A_M}{B_M} \stackrel{d}{=} X, 
\end{equation}
where $N \sim U([\![1,M]\!])$ and $X$ is a random variable with a non-degenerate distribution, i.e. $X$ attains at least two different values with positive probability. By Proposition \ref{LemRepresentationf} and since $B_M \to \infty$, we can assume that $f$ is of the form

\[f(x) = \left ( \iota(x) - \frac{1}{2} \right) \sum_{i=1}^{\nu} H_i + \sum_{i=1}^{\nu}H_i \left(  \mathbb{1}_{[0, \gamma_i)}(x) -\iota(\gamma_i)  \right)\]

where $H_i \in \R$. Let $(n_\ell)_{\ell \in \N}$ be the sequence of integers from Lemma \ref{indicator_lem} and, for some $c \in (0,1]$, define $M_{\ell} := \lfloor ca_{n_{\ell} +1} \rfloor q_{n_{\ell}} + q_{\kj}-1$.
Clearly, any $N \in [0,M_{\ell}]$ has a unique representation of the form $N = b_{\ell}q_{n_{\ell}} +N' $
where $0 \leq b_{\ell} \leq \lfloor ca_{n_{\ell}+1}\rfloor $ and $0 \leq N' \leq q_{n_{\ell}}-1$.
It follows immediately from the definition that we can decompose the Birkhoff sum as

\[
\begin{split}S_N(f,\alpha,x_0) &= S_{b_{\ell}q_{\ell}}(f,\alpha,x_0) +S_{N'}\left(f, \alpha, T_{\alpha}^{b_\ell  q_{n_{\ell}}}(x_0)\right)
\\&= S_{b_{\ell}q_{\ell}}(f,\alpha,x_0) + S_{N'}(f, \alpha, x_0 +  b_\ell  q_{n_{\ell}} \alpha)
. \end{split}\]

Applying the Denjoy-Koksma inequality (see \eqref{denj_koks}) shows that

\begin{equation*}\lvert S_{N'}(f, \alpha, x_0 +  b_\ell  q_{n_{\ell}} \alpha) \rvert
\ll_f \sum_{i = 1}^{n_{\ell}} a_i,\end{equation*}

which by the properties of $(\kj)_{\ell \in \N}$ implies that
\begin{equation*}
\frac{S_{N'}(f, \alpha, x_0 +  b_\ell  q_{n_{\ell}} \alpha)   }{a_{\kj +1}} = o(1), \quad \ell \to \infty.
\end{equation*}

If $N_{\ell} \sim U([\![0, M_{\ell}]\!])$, then it is easy to see that
\begin{equation*}
N_{\ell } \stackrel{d}{=} b_{\ell}q_{n_{\ell}} +N',
\end{equation*}
where $b_{\ell} \sim U([\![0,\lfloor ca_{n_{\ell}+1}\rfloor]\!])$, $N' \sim U([\![0,q_{n_{\ell}}-1]\!])$ and $b_{\ell}$ and $N'$ are independent. Hence,
\begin{equation*}
    \frac{S_{N_\ell}( f, \alpha, x_0)}{a_{\kj +1}} \stackrel{d}{=}  \frac{S_{b_\ell q_{\kj }}( f, \alpha, x_0)}{a_{\kj +1}} + o(1).
\end{equation*}
Thus we get for any $x \in \R$
\begin{align*}
\frac{1}{M_{\ell}} \# \left \{  1 \leq N \leq M_{\ell} : \frac{S_N(f, \alpha, x_0)}{a_{n_{\ell}+1}} \leq x \right \} 
&= \frac{1}{M_{\ell}} \# \left \{  0 \leq N \leq M_{\ell} : \frac{S_N(f, \alpha, x_0)}{a_{n_{\ell}+1}} \leq x \right \} + o(1)
\\& = \mathbb{P} \left[ \frac{S_{N_\ell}(f, \alpha, x_0)}{a_{n_\ell +1}} \leq x \right] + o(1) \\
&= \mathbb{P} \left[ \frac{S_{b_\ell q_{\kj }}(f, \alpha, x_0)}{a_{n_\ell +1}} \leq x + o(1)\right] + o(1) \\
&= \mathbb{P} \left[ \frac{S_{ \lfloor U_c a_{\kj +1} \rfloor q_{\kj }}(f, \alpha, x_0)}{a_{n_\ell +1}} \leq x + o(1)\right] + o(1),
\end{align*}
where $U_c \sim U([0,c])$. In the last line, we used that
\[
\mathbb{P} \left[ \frac{S_{b_\ell q_{\kj }}(f, \alpha, x_0)}{a_{n_\ell +1}} \leq y \right] = \mathbb{P} \left[ \frac{S_{\lfloor U_c a_{\kj +1} \rfloor q_{\kj }}(f, \alpha, x_0)}{a_{n_\ell +1}} \leq y \right]  + o(1),
\]
uniformly in $y \in \R$. Moreover, by Lemma \ref{indicator_lem} we get the (almost sure) limit
\begin{align*}
 \lim_{\ell \to \infty} \frac{S_{ \lfloor U_c a_{n_{\ell}+1} \rfloor q_{n_{\ell}} }(f, \alpha, x_0) }{a_{n_\ell +1}} &= g(U_c),
 \end{align*}

 where, for $x \in [0,1]$, \[g(x) := 
 \left(\sum_{i = 1}^{\nu} H_i\right)\left(\int_{0}^{x} \iota( y + \overline{x_0}) \, \mathrm{d}y - \frac{x}{2}\right) + \sum_{i = 1}^{\nu} H_i\left(\int_{0}^{x} \mathds{1}_{\left[0, \overline{\gamma_i}\right]}\left(y + \overline{x_0}\right) \,\mathrm{d}y - x\iota(\overline{\gamma_i})\right).\]

Since $g(U_c)$ has a continuous distribution, this implies that
\[
\lim_{\ell \rightarrow \infty} \frac{1}{M_{\ell}} \# \left \{  1 \leq N \leq M_{\ell} : \frac{S_N(f, \alpha, x_0)}{a_{n_{\ell}+1}} \leq x \right \} = \mathbb{P} \left [ g(U_c) \leq x \right ] .
\]
Now let $ \tilde{A}_M := 0$ and $\tilde{B}_M := \frac{M}{q_{n(M)}}$, where 
$q_{n(M)} \leq  M < q_{n(M) +1}$.
We have shown in the previous argument that, for any $c \in (0,1]$ and for $ (n_\ell)_{\ell \in \N}$ as before, we have
\begin{equation*}
     \lim_{\ell \to \infty} \frac{S_{ \lfloor U_c a_{n_{\ell}+1} \rfloor q_{n_{\ell}} }(f, \alpha, x_0) - \tilde{A}_{M_\ell} }{\tilde{B}_{M_\ell}} \stackrel{d}{=} c g(U_c).
\end{equation*}
By the convergence of types theorem (see, e.g., \cite[Theorem 14.2]{billingsley}) and since the limit in \eqref{LiminDist} also holds along every subsequence tending to infinity, there exist quantities $B_c > 0$ and $A_c \in \R$ such that for any $c \in (0,1]$ we have
\begin{equation*}
    c g(U_c) \stackrel{d}{=} B_c X + A_c.
\end{equation*}
This implies that for any $0 < c_1, c_2 \leq 1$, we can write 
\begin{equation}
\label{EqEqinDist}
g(U_{c_1}) \stackrel{d}{=} B(c_1, c_2) g(U_{c_2}) + A(c_1, c_2),
\end{equation}
where $B(c_1, c_2) > 0$ and $A(c_1, c_2) \in \R$.

We now collect a few properties of the function $g(x)$ for $x \in \R$. First, we note that $g(0)=g(1) = 0$. Further, $g$ is differentiable except in all points of the form $ \overline{\gamma}_i+ \overline{x}_0$ and $g$ is non-constant. To see the latter, we fix $ \delta > 0$ small enough such that $ \delta < \min_{i=2, \ldots, \nu} \lVert  \overline{\gamma}_1- \overline{\gamma_i} \rVert$ (which is possible because $\overline{\gamma_1} \neq \overline{\gamma_i}$ for all $i = 2, \ldots, \nu$). We then get
\begin{equation*}
g'\left ( \iota(\overline{\gamma}_1 + \overline{x}_0) - \frac{\delta}{2} \right) - g'\left ( \iota(\overline{\gamma}_1 + \overline{x}_0) + \frac{\delta}{2} \right) = \delta \left( \sum_{i=1}^{\nu} H_i \right) +H_1.
\end{equation*}
By choice of $f$, there exists at least one $H_i \neq 0$, thus, we may assume $H_1 \neq 0$. Since $\delta$ can be chosen arbitrarily small, it follows that $g'$ is not constant and hence $g$ is not constant. Hence, locally to the right of $0$, $g(x)$ is either monotonically increasing or monotonically decreasing. In the following we discuss the case where $g(x)$ is increasing, the case where $g(x)$ is decreasing can be handled analogously. It follows that there exists an $\varepsilon \in (0,1)$ and $ \delta \in (0,1]$ with $ \varepsilon < \delta$ with the following properties: The function $g$ is increasing on $[0, \epsilon]$ with $g ( \varepsilon) > 0$. On $[\varepsilon, \delta]$, $g$ is decreasing and $ 0 \leq g( \delta) < g( \varepsilon)$. 
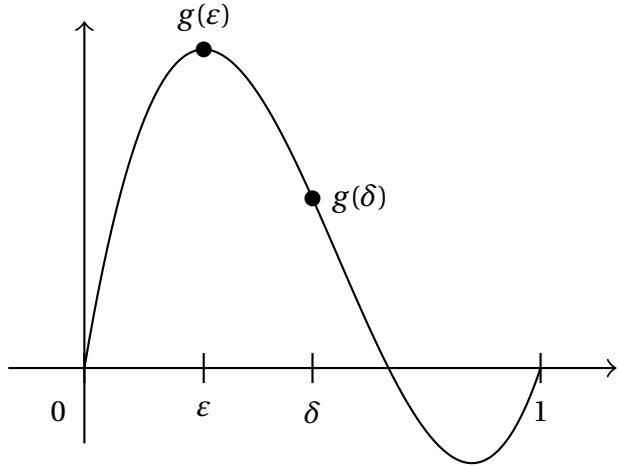
\begin{figure}[H]
        \centering
        \scalebox{2}{
            \begin{tikzpicture}
                \draw[->] (0,-0.5) -- (0,2.3);
                \draw[->] (-0.5,0) -- (3.5,0);
                \draw[domain=0:3] plot[samples=100] (\x,{\x*\x*\x-5*\x*\x+6*\x});
                \fill[] (0.78475,2.11261) circle[radius=1.5pt] node[above] {\tiny $ g(\varepsilon) $};
                \fill[] (1.5,1.125) circle[radius=1.5pt] node[right] {\tiny $ g(\delta) $};
                \draw[] (0,0.1) -- (0,-0.1) node[below left] {\tiny $ 0 $};
                \draw[] (0.78475,0.1) -- (0.78475,-0.1) node[below] {\tiny $ \varepsilon $};
                \draw[] (1.5,0.1) -- (1.5,-0.1) node[below] {\tiny $ \delta $};
                \draw[] (3,0.1) -- (3,-0.1) node[below] {\tiny $ 1 $};
            \end{tikzpicture}
        }
        \captionsetup{labelformat=empty}
        \caption{Illustration of the argument above. Clearly, $g([0,\varepsilon]) = g([0,\delta])$.}
        \label{fig:einplot}
    \end{figure} 

Using \eqref{EqEqinDist} we infer
\begin{equation*}
    g(U_{\varepsilon}) \stackrel{d}{=} B(\varepsilon, \delta)  g(U_{\delta}) + A(\varepsilon, \delta).
\end{equation*}

However, by the choice of $\varepsilon$ and $\delta$, we have $g([0, \varepsilon]) = g([0, \delta])$, which immediately implies that $A(\varepsilon, \delta)=0$ and $B(\varepsilon, \delta)=1$.
By construction, we have 

\[\mathbb{P}\left[g(U_{\varepsilon}) > g(\delta)\right] < \mathbb{P}\left[g(U_{\delta}) > g(\delta)\right],\]

which is an immediate contradiction to $g(U_\varepsilon) \stackrel{d}{=} g(U_\delta)$.
\end{proof}

\subsection*{Acknowledgements} 
We would like to thank Bence Borda for many valuable discussions.
LF and MH are supported by the Austrian Science Fund (FWF) Project P 35322 \textit{Zufall und Determinismus in Analysis und Zahlentheorie}.


\begin{thebibliography}{HD82}


\bibitem{aake}
J. Aaronson, M. Keane, \textit{The visitors to zero of some deterministic random walks},  Proc. London Math. Soc. Third Series, s3-44(3), p. 535--553, 1982.

\bibitem{all_shall}
J. Allouche, J. Shallit,
{\it Automatic Sequences: Theory, Applications, Generalizations}, 
Cambridge University Press, 2003.

\bibitem{addo}
A. Avila, D. Dolgopyat, E. Duryev, O. Sarig, \textit{The visits to zero of a random walk driven by an irrational rotation}, Isr. J. Math. 207, p. 653--717, 2015.

\bibitem{beck1}
J. Beck, {\it Probabilistic Diophantine approximation: Randomness in lattice point counting}, Springer Monographs in Mathematics, Springer, Cham, International Publishing, 2014.

\bibitem{beck2}
J. Beck, {\it Randomness of the square root of 2 and the Giant Leap, part 1}, Period. Math. Hungar. 60, p. 137--242, 2010. https://doi.org/10.1007/s10998-010-2137-9

\bibitem{beck3}
J. Beck, {\it Randomness of the square root of 2 and the giant leap, part 2.} Period. Math. Hungar. 62, p. 127--246, 2011. https://doi.org/10.1007/s10998-011-6127-3


\bibitem{borda}
B. Borda, {\it On the distribution of Sudler products and Birkhoff sums for the irrational rotation}, pre-print: arXiv:2104.06716, to appear in: Ann. Inst. Fourier (Grenoble) 


\bibitem{brom_ulc}
M. Bromberg, C. Ulcigrai, {\it A temporal central limit theorem for real-valued cocycles over rotations}, Ann. inst. Henri Poincare (B) Probab. Stat., 54(4), p. 2304--2334, 2018. https://doi.org/10.1214/17-aihp872


\bibitem{diamond_vaaler}
H. Diamond, J. Vaaler,
{\it Estimates for partial sums of continued fraction partial quotients},
Pacific J. Math. 122, p. 73--82, 1986.

\bibitem{dol_sar}
D. Dolgopyat, O. Sarig, \textit{Temporal distributional limit theorems for dynamical systems}, J. Stat. Phys., 166, p. 680--713, 2017. https://doi.org/10.1007/s10955-016-1689-3

\bibitem{dol_sar_no}
D. Dolgopyat, O. Sarig, \textit{No temporal distributional limit theorem for a.e. irrational translation},
Ann. H. Lebesgue 1, p. 127--148, 2018.

\bibitem{quenched}
D. Dolgpyat, O. Sarig, \textit{Quenched and annealed temporal limit theorems for circle rotations}, Asterisque, 415, p. 59--85, 2020. https://doi.org/10.24033/ast.11100

\bibitem{dol_fay}
D. Dolgopyat, B. Fayad, \textit{Limit theorems for toral translations}, Proc. Symp. Pure Math., 89, p. 227--277, 2015. 

\bibitem{duff_sch}
R. Duffin, A. Schaeffer,  \textit{Khintchine’s problem in metric Diophantine approximation}, Duke Math. J., 8(2), p. 243--255, 1941. https://doi.org/10.1215/s0012-7094-41-00818-9


\bibitem{upper_dens}
L. Fruehwirth, M. Hauke, \textit{On the metric upper density of Birkhoff sums for irrational rotations}, preprint:arXiv:2303.15992


\bibitem{herman}
M.R. Herman, \textit{Sur la Conjugaison Différentiable des Difféomorphismes du Cercle a des Rotations}, Inst. Hautes Études Sci. Publ. Math. No. 49, 5--233, 1979. https://doi.org/10.1007/bf02684798


\bibitem{Kesten1960}
H. Kesten,
{ \it Uniform distribution mod 1},
Ann. of Math., p. 445--471, 1960.


\bibitem{kuipers}
L. Kuipers, H. Niederreiter, {\it Uniform Distribution of Sequences}, Wiley, 1974.


\bibitem{schmidt}
K. Schmidt. \textit{A cylinder flow arising from irregularity of distribution},
Compos. Math. 36(3), p. 225--232, 1978.


\bibitem{rock_sz}
 A. M. Rockett, P. Sz\"usz,
 {\it Continued fractions}, World Scientific Publishing, River Edge, NJ, 1992.

\bibitem{gall}
P. Gallagher, {\it Approximation by reduced fractions}, J. Math. Soc. Japan, p. 342--345, 1961. https://doi.org/10.2969/jmsj/01340342

\bibitem{billingsley}
P. Billingsley,{ \it Probability and measure}, John Wiley \& Sons, Inc., New York, 1995.
	    
\end{thebibliography}
\end{document}